\documentclass{dcds-sEA} 
\usepackage{amsmath}
\usepackage{amssymb}
\usepackage{paralist}
\usepackage{stfloats}
\usepackage[misc]{ifsym}
\usepackage{comment}
\usepackage{graphicx}%
\usepackage{subfigure}
\usepackage{multirow}%
\usepackage{amsfonts}%
\usepackage{amsthm}
\usepackage{mathtools}
\usepackage{cases}
\usepackage{mathrsfs}%
\usepackage[title]{appendix}%
\usepackage{xcolor}%
\usepackage{textcomp}%
\usepackage{manyfoot}%
\usepackage{booktabs}%
\usepackage{algorithm}%
\usepackage{algorithmicx}%
\usepackage{algpseudocode}%
\usepackage{listings}%
\usepackage{latexsym}
\usepackage{indentfirst}
\usepackage{tabularx}
\usepackage{multirow} 
\usepackage{adjustbox}

\usepackage{epsfig} 
\usepackage{epstopdf} 
\usepackage[colorlinks=true]{hyperref}
\hypersetup{urlcolor=blue, citecolor=red}

\allowdisplaybreaks

\def\currentvolume{X}
 \def\currentissue{X}
  \def\currentyear{200X}
   \def\currentmonth{XX}
    \def\ppages{X-XX}
     \def\DOI{10.3934/dcdss.2024202}

\newtheorem{theorem}{Theorem}[section]

\newtheorem{proposition}[theorem]{Proposition}

\theoremstyle{definition}

\title[Symmetry classification and invariant solutions of the GMR process]
{Symmetry classification and invariant \\solutions of the classical geometric \\mean reversion process} 

\author[Jin Zhang and Dapeng Gao]{}


\thanks{The first author is supported by [EDJP grant JJKH20211030KJ]}

\thanks{$^*$Corresponding author: Jin Zhang}

\begin{document}
\maketitle

\centerline{\scshape
Jin Zhang$^{{\href{mailto:jinzhang@jlu.edu.cn}{\textrm{\Letter}*}}}$
and Dapeng Gao$^{{\href{mailto:gaodp22@mails.jlu.edu.cn}{\textrm{\Letter}}}}$}

\medskip

{\footnotesize
 \centerline{School of Mathematics, Jilin University, Changchun, Jilin 130012, China}
} 



\bigskip



\begin{abstract}
Based on the Lie symmetry method, we investigate a Feynman-Kac formula for the classical geometric mean reversion process, which effectively describing the dynamics of short-term interest rates. The Lie algebra of infinitesimal symmetries and the corresponding one-parameter symmetry groups of the equation are obtained. An optimal system of invariant solutions are constructed by a derived optimal system of one-dimensional subalgebras. Because of taking into account a supply response to price rises, this equation provides for a more realistic assumption than the geometric Brownian motion in many investment scenarios.
\end{abstract}



\section{Introduction}
\label{sec:sample1}

Over the past few decades, stochastic volatility models have gained popularity in the domain of derivative pricing and hedging. Empirical evidence indicates that asset volatility is stochastic, and the traditional Black-Scholes-Merton (BSM) model \cite{black1973pricing} falls short in capturing important aspects of asset returns, such as skewness, leptokurtosis, and conditional heteroscedasticity.

\medskip
Notable contributions in stochastic volatility modeling include the works of Hull and White \cite{hull1987pricing}, Scott \cite{scott1987option}, and Wiggins \cite{wiggins1987option}, who employed numerical techniques to solve two-dimensional partial differential equations for pricing solutions. Stein and Stein \cite{stein1991stock} were able to analytically value options by assuming that volatility is uncorrelated with the spot asset. However, it is known that volatility shocks are negatively correlated with asset price shocks, resulting in a leverage effect where an increase in volatility leads to a decrease in stock prices, and vice versa. This correlation is crucial for capturing the skewed distribution of asset prices, and without it, pricing models cannot fully account for such skewness effects.

\medskip
A notable stochastic volatility model in finance is the mean-reverting model proposed by Heston \cite{heston1993closed}, which provides an exact solution for pricing European options. There has been a recent surge in research on the use of symmetry analysis for option pricing differential equations (see \cite{bordag2007explicit, bordag2011models, caister2010solving, caister2011optimal, craddock2009calculation, wang2014optimal, yang2007symmetry} and the references therein for example). Additionally, there are further investigations on optimal systems and corresponding group invariant solutions for different partial differential equations related to finance (see \cite{chou2001note, chou2004optimal,kaibe2021application, lo2001valuation, lo2006lie} for example).

\medskip
In the context of mean reversion, we first consider the classical Ornstein-Uhlenbeck (OU) process, which has the disadvantage of potentially taking negative values with positive probability, posing challenges for financial applications requiring non-negativity. Khalique and Motsepa \cite{khalique2018lie} conducted a symmetry analysis of the Vasicek pricing equation, which shares similar limitations to the OU process despite its frequent use in interest rate theory.

\medskip
In contrast, the Cox-Ingersoll-Ross (CIR) model \cite{cox85} addresses this issue by ensuring non-negativity while effectively capturing the stochastic nature of interest rates. In \cite{sinkala2008optimal, sinkala2008zero}, the authors analyzed the valuation partial differential equations for European interest rate derivatives within the CIR framework, deriving group-invariant solutions for zero-coupon bond pricing. In \cite{tang2022lie}, the authors exploited Lie symmetry analysis to solve the generalized CIR model.

\medskip
As is well known, in the case of commodities, interest rates and exchange rates, a geometric mean reversion (GMR) model holds more economic rationale compared to the geometric Brownian motion (GBM). The GMR process effectively incorporates both mean-reverting dynamics and stochastic volatility, making it a compelling framework for further investigation.

\medskip
Thus, we will focus on the GMR process, which retains the advantages of non-negativity and provides a more realistic representation of asset behavior.
The classical GMR process takes the form:
\begin{equation}\label{sde}
{\rm d}X_t=k(\alpha-X_t)X_t{\rm d}t+\sigma X_t{\rm d}W_t, \quad X_0=x,
\end{equation}
where $X_t$ represents the stock prices, commodity prices, interest rates, etc., $\alpha\geqslant0$ is the long-run equilibrium level of $X_t$, $k>0$ represents the speed of reversion, $\sigma>0$ is the volatility in the price process, and $W_t$ is the standard Brownian motion \cite{karatzas2012brownian} (also known as the Wiener process \cite{ahmad1988introduction, gardiner1985handbook}).

\medskip
Merton \cite{merton1975asymptotic} investigated a growth model with uncertainty and considered the special case of a Cobb-Douglas production function with a constant savings function. The process (\ref{sde}) naturally represented the output-to-capital ratio in this context. Additionally, this process has found application in classical investment theory and real options as a realistic model for the value of project with mean-reverting returns.

\medskip
Dixit and Pindyck \cite{dixit1994investment} offered a thorough discussion of investment model where the value function followed the GMR process. However, they didn't provide a rationale for why the value function was mean-reverting.

\medskip
Metcalf and Hassett \cite{metcalf1995investment} proposed a different discussion of the GMR model by modeling mean reversion in the output price and providing a rationale for why the price mean-reverts. In addition, they compared the investment strategies based on the GBM and the mean reversion. Despite its greater complexity, the mean reversion was considered a more realistic assumption in most cases as it allowed for supply responses to increasing prices. They demonstrated that the cumulative investment was generally unaffected when employing a mean reversion process instead of the GBM, and they explained this outcome.

\medskip
Ewald and Yang \cite{ewald08} applied this process to model the discounted project value. The mean reverting structure of the project value process enriched the model and had greater economic significance. Additionally, this process is relevant in other aspects of finance such as interest rate theory, stochastic volatility or option pricing, see for example \cite{de2015modeling, dung2011fractional, ewald2006new, yang2010} and references therein.

\medskip
The Feynman-Kac formula \cite{deck2002white, feynman1948space, kac1949distributions, moral2004feynman} establishes a link between the realm of partial differential equations and the pertinent probability theory associated with stochastic differential equations. A notable application of this formula is exemplified in the derivation of the Black-Scholes-Merton model,
in which the partial differential equation is transformed into the heat equation via variable substitutions.

\medskip
For the classical GMR process characterized by (\ref{sde}), we utilize the Feynman-Kac formula (see \cite{durrett1984brownian,durrett1985} for example) and an appropriate function $c(\cdot)$ to obtain the value formula:
\begin{equation*}
    u(t,x) = E_x\left[e^{\int_0^t c(X_s) {\rm d}s}\right],
\end{equation*}
where $X_t$ is the GMR process given by (\ref{sde}) with $X_0=x$, and $u$ satisfies the parabolic partial differential equation \cite{oksendal2003,shreve2004stochastic}:
\begin{eqnarray}\label{eq2}
    u_t = \frac{1}{2}\sigma^2 x^2u_{xx} + kx(\alpha-x)u_x + c(x)u.
\end{eqnarray}

\medskip
In this context, the choice of $c(x)=\lambda x^2$ introduces a term that reflects both volatility and mean reversion behavior in the model. This form captures the fact that volatility often increases with the magnitude of $x$, a common observation in financial markets, where larger asset values tend to exhibit greater fluctuations. Moreover, the quadratic term $\lambda x^2$ can also represent a form of reinforcement or decay in the system: for $\lambda>0$, it models growth, while for $\lambda<0$, it introduces a damping effect that encourages the system to stabilize around the equilibrium. Such dynamics are particularly useful in scenarios where asset prices are expected to revert to a long-term mean, such as in commodities or interest rate models. Thus, our focus will be on the following partial differential equation:
\begin{eqnarray}\label{eq1}
    u_t = \frac{1}{2}\sigma^2 x^2u_{xx} + kx(\alpha-x)u_x + \lambda x^2u,
\end{eqnarray}
which serves as the governing differential equation for the GMR process.

\medskip
In this paper, we investigate the existence of symmetry group for equation (\ref{eq1}) by using symmetry analysis (see \cite{bluman02, bluman89, ibragimov1999elementary, olver93, olver95,ovsiannikov2014group} for example). The symmetry group of a differential equation is a set of transformations that map solutions of the equation into other solutions. Once the symmetry group is determined, various applications become possible. One can utilize the defining property of the group to construct new solutions from known ones and obtain corresponding invariant solutions by symmetry reductions. Since there are usually infinitely many invariant solutions, an optimal system is required for categorize all potential invariant solutions for equation (\ref{eq1}).

\medskip
The framework of this article is structured as follows: In Section \ref{sec2}, we derive the Lie algebra of infinitesimal symmetries and the corresponding one-parameter symmetry groups for (\ref{eq1}). In Section \ref{sec3}, an optimal system for the symmetry algebra is constructed via the commutation relations and adjoint representation of subalgebras, leading to intriguing exact solutions. In Section \ref{sec4}, symmetry reductions are considered and novel explicit solutions of (\ref{eq1}) are obtained by employing the one-dimensional optimal system of the symmetry algebra.

\section{Symmetry analysis} \label{sec2}

In this section, we apply Lie's classical method to determine the symmetry groups of (\ref{eq1}). We begin by considering a one-parameter group of infinitesimal transformations:
\begin{equation}\label{trans}
    \begin{aligned}
     \tilde t\,=\,t+\varepsilon\tau(t,x,u)+O(\varepsilon^2),\\
     \tilde x=x+\varepsilon\xi(t,x,u)+O(\varepsilon^2),\\
     \tilde u=u+\varepsilon\phi(t,x,u)+O(\varepsilon^2),
    \end{aligned}
\end{equation}
with a continuous parameter $\varepsilon$. The infinitesimal generator associated with group of transformations (\ref{trans}) can be expressed as:
\begin{equation*}
    V=\tau(t,x,u)\partial_t+\xi(t,x,u)\partial_x+\phi(t,x,u)\partial_u.
\end{equation*}
The requirement that (\ref{eq1}) is invariant under the transformations (\ref{trans}) is determined by the infinitesimal criterion of invariance. This criterion states that the action of the prolongation of the infinitesimal generator on the equation must be identically zero, modulo the equation under consideration. Here, we require that
\begin{eqnarray}\label{0}
    pr^{(2)}V\left(\frac{1}{2}\sigma^2x^2u_{xx}+kx(\alpha-x)u_x+\lambda x^2u-u_t\right)\Big|_{(\ref{eq1})}=0,
\end{eqnarray}
where $pr^{(2)}V$ is the second prolongation of vector field $V$, which has the form
\begin{eqnarray*}
    pr^{(2)}V=V+\phi^t\partial_{u_t}+\phi^x\partial_{u_x}+\phi^{tt}\partial_{u_{tt}}+\phi^{tx}\partial_{u_{tx}}+\phi^{xx}\partial_{u_{xx}},
\end{eqnarray*}
where $\phi^J=D_J(\phi-\tau u_t-\xi u_x)+\tau u_{J,t}+\xi u_{J,x}$, $J$ represents the multi-indices of orders $1\leq \#J\leq2$.
Substituting the general formulas of $\phi^J$ into (\ref{0}) and equating the coefficients of the various monomials for the partial derivatives of $u$, we immediately obtain an overdetermined set of equations, which are called the determining equations
\begin{eqnarray*}
    & & \tau_x=0, \ \ \tau_u=0, \ \ \xi_u=0, \ \ \phi_{uu}=0, \ \ 2\xi+x(\tau_t-2\xi_x)=0, \label{1} \\
    & & k(\alpha-2x)\xi+\xi_t+kx(\alpha-x)(\tau_t-\xi_x)+\frac{1}{2}\sigma^2x^2(2\phi_{xu}-\xi_{xx})=0, \label{2} \\
    & & 2\lambda xu\xi-\phi_t+kx(\alpha-x)\phi_x+\frac{1}{2}\sigma^2x^2\phi_{xx}+\lambda x^2\phi-\lambda x^2u(\phi_u-\tau_t)=0, \label{3}
\end{eqnarray*}
where the subscripts denote differentiation with respect to the indicated variables. Solving the determining equations, we then obtain the coefficient functions
\begin{equation*}
	\begin{aligned}
		\tau=\ &c_1+c_2t+c_3t^2, \\
		\xi=\ &\frac{c_2}{2}x\ln x+c_3tx\ln x+c_4x+c_5tx, \\
		\phi=\ &c_2 \! \left(\! \frac{k}{2\sigma^2}x\ln x \!+\! \frac{1}{4}\ln x \!-\! \frac{\sigma^2}{8}t \! \right)\! u \!+\! c_3 \! \left(\! \frac{k}{\sigma^2}tx\ln x \!-\! \frac{1}{2\sigma^2}\ln^2x \!+\! \frac{1}{2}t\ln x \!-\! \frac{\sigma^2}{8}t^2 \!-\! \frac{1}{2}t \! \right)\! u \\
		&\!+c_4\frac{k}{\sigma^2}xu+c_5\left(\frac{k}{\sigma^2}tx-\frac{1}{\sigma^2}\ln x+\frac{1}{2}t\right)u+c_6u+\varphi(t,x),
	\end{aligned}
\end{equation*}
where $c_1$, $\cdots$, $c_6$ are arbitrary constants for $\alpha=0$ and $\lambda=\frac{k^2}{2\sigma^2}$, and $c_2=c_3=c_4=c_5=0$, $c_1$, $c_6$ are arbitrary constants for other cases.

\medskip
Summarizing the above discussions, we can immediately derive two results concerning the symmetry algebra of (\ref{eq1}).
\begin{proposition}
For the case  $\alpha=0$ and $\lambda=\frac{k^2}{2\sigma^2}$, the Lie algebra of infinitesimal symmetries of (\ref{eq1}) is spanned by the six vector fields
\begin{flalign*}
		&V_1=\partial_t, \\
		&V_2=t\partial_t+\frac{1}{2}x\ln x\partial_x+\left(\frac{k}{2\sigma^2}x\ln x
		+\frac{1}{4}\ln x-\frac{\sigma^2}{8}t\right)u\partial_u, \\
		&V_3=t^2\partial_t+tx\ln x\partial_x+\left(\frac{k}{\sigma^2}tx\ln x-\frac{1}{2\sigma^2}\ln^2x
		+\frac{1}{2}t\ln x-\frac{\sigma^2}{8}t^2-\frac{1}{2}t\right)u\partial_u, \\
		&V_4=x\partial_x+\frac{k}{\sigma^2}xu\partial_u, \\ &V_5=tx\partial_x+\left(\frac{k}{\sigma^2}tx
		-\frac{1}{\sigma^2}\ln x+\frac{1}{2}t\right)u\partial_u, \\
		&V_6=u\partial_u
\end{flalign*}
and the infinite-dimensional subalgebra $V_{\varphi}=\varphi(t,x)\partial_u$, where $\varphi(t,x)$ is an arbitrary solution of (\ref{eq1}).
\end{proposition}
\begin{proposition}
	For the case $\alpha\neq0$ or $\lambda\neq\frac{k^2}{2\sigma^2}$, the Lie algebra of infinitesimal symmetries of (\ref{eq1}) is spanned by the two vector fields
	$$V_1=\partial_t,~~~~~V_6=u\partial_u$$
	and the infinite-dimensional subalgebra $V_{\varphi}=\varphi(t,x)\partial_u$, where $\varphi(t,x)$ is an arbitrary solution of (\ref{eq1}).
\end{proposition}

The corresponding one-parameter symmetry groups generated by $V_i$ are then:
\begin{align*}
    &\ G_1: \left(t, x, u\right)\rightarrow\left(t+\varepsilon, x, u\right),\\
    &\ G_2: \left(t, x, u\right)\rightarrow\left(e^\varepsilon t, e^{e^{\frac{\varepsilon}{2}}\ln x}, ue^{\frac{k}{\sigma^2}\big(e^{e^{\frac{\varepsilon}{2}}\ln x}-x\big)+\frac{1}{2}(e^{\frac{\varepsilon}{2}}-1)\ln x-\frac{\sigma^2}{8}t(e^{\varepsilon}-1)}\right),\\
    &\ G_3: \left(t, x, u\right)\rightarrow\left(\frac{t}{1-\varepsilon t}, e^{\frac{\ln x}{1-\varepsilon t}}, ue^{\frac{k}{\sigma^2}\big(e^{\frac{\ln x}{1-\varepsilon t}}-x\big)+\frac{\varepsilon t}{1-\varepsilon t}\left( \frac{\ln^2x}{2\sigma^2t}-\frac{1}{2}\ln x+\frac{\sigma^2}{8}t \right)+\frac{1}{2}\ln( 1-\varepsilon t)}\right),\\
    &\ G_4: \left(t, x, u\right)\rightarrow\left(t, e^\varepsilon x, ue^{\frac{k}{\sigma^2}x\left( e^\varepsilon -1 \right)} \right),\\
    &\ G_5: \left(t, x, u\right)\rightarrow\left(t, e^{\varepsilon t} x, ue^{\frac{k}{\sigma^2}x( e^{\varepsilon t}-1)-\frac{\varepsilon^2}{2\sigma^2}t-\frac{\varepsilon}{\sigma^2}\ln x+\frac{\varepsilon t}{2}} \right),\\
    &\ G_6: \left(t, x, u\right)\rightarrow\left(t, x, e^\varepsilon u\right),\\
    &\ G_{\varphi}: \left(t, x, u\right)\rightarrow\left(t, x, u+\varepsilon\varphi(t,x)\right).
\end{align*}

\medskip
The group $G_1$ represents the time translation invariance of (\ref{eq1}). The symmetry groups $G_6$ and $G_{\varphi}$ reflect the linearity of (\ref{eq1}), indicating that the equation's solution remains closed under addition and multiplication by a constant. The symmetry group in the case of Proposition 2.2 is trivial, hence we only consider the case of Propositon 2.1 in the remaining discussions.

\section{Optimal system of one-dimensional symmetry algebras} \label{sec3}

Olver \cite{olver93} and Ovsiannikov \cite{ovsiannikov2014group} proposed effective methods for categorizing the subalgebras of the Lie algebra formed by symmetries from different viewpoints in order to classify all invariant solutions. In this section, we will construct an optimal system for the one-dimensional subalgebras of (\ref{eq1}).

\medskip
Since there are no invariant solutions for the infinite-dimensional subalgebra $V_{\varphi}$, as Olver \cite{olver93} explains, the classification problem will not take it into account. Moving on, we shall consider the Lie algebra $\mathfrak{g}$, which is spanned by $\{V_1, V_2, V_3, V_4, V_5, V_6\} $.

\medskip
To classify invariant solutions, we must initially examine the commutation relations and the adjoint representation of the Lie algebra $\mathfrak{g}$ of infinitesimal symmetries of (\ref{eq1}), which yields the optimal system.

\medskip
The commutation relations between vector fields $V_i$ is given in Table \ref{tab1} with the entry in row $i$ and column $j$ representing the Lie bracket $[V_i,V_j]=V_iV_j-V_jV_i$.
\begin{table}[H]
\,\caption{Commutator table of Lie algebra for (\ref{eq1}).}
\begin{center}
	\begin{tabular}{c|cccccc}
		\hline
		$[V_i,V_j]$ & $V_1$ & $V_2$ & $V_3$ & $V_4$ & $V_5$ & $V_6$ \\ \hline
		$V_1$ & 0 & $V_1-\frac{\sigma^2}{8}V_6$ & $2V_2-\frac{1}{2}V_6$ & 0 & $V_4+\frac{1}{2}V_6$ & 0 \\
		$V_2$ & $-V_1+\frac{\sigma^2}{8}V_6$ & 0 & $V_3$ & $-\frac{1}{2}V_4-\frac{1}{4}V_6$ & $\frac{1}{2}V_5$ & 0 \\
		$V_3$ & $-2V_2+\frac{1}{2}V_6$ & $-V_3$ & 0 & $-V_5$ & 0 & 0 \\
		$V_4$ & 0 & $\frac{1}{2}V_4+\frac{1}{4}V_6$ & $V_5$ & 0 & $-\frac{1}{\sigma^2}V_6$ & 0 \\
		$V_5$ & $-V_4-\frac{1}{2}V_6$ & $-\frac{1}{2}V_5$ & 0 & $\frac{1}{\sigma^2}V_6$ & 0 & 0 \\
		$V_6$ & 0 & 0 & 0 & 0 & 0 & 0 \\
		\hline
	\end{tabular}
\end{center}
\label{tab1}
\end{table}

\medskip
The adjoint representation can be expressed in the form of Lie series
\begin{equation}\label{14}
    {\rm Ad}(\exp(\varepsilon V_i))V_j=V_j-\varepsilon [V_i, V_j]+\frac{\varepsilon^2}{2}[V_i,[V_i,V_j]]-\cdots.
\end{equation}
By assisting with Table \ref{tab1} and the Lie series (\ref{14}), the adjoint representation table of (\ref{eq1}) can be constructed in Table 2 with the $(i,j)$-th entry indicating ${\rm Ad}(\exp(\varepsilon V_i))V_j$.
\begin{table}[H]
\,\caption{Adjoint representation of subalgebras.}
\resizebox{\textwidth}{11.5mm}{
\centering
    \begin{tabular}{c|ccccccc}
 \hline
	Ad & $V_1$ & $V_2$ & $V_3$ & $V_4$ & $V_5$ & $V_6$  \\ \hline
        $V_1$ & $V_1$ & $V_2-\varepsilon(V_1-\frac{\sigma^2}{8}V_6)$ & $V_3+\varepsilon^2 V_1-2\varepsilon V_2+\frac{4\varepsilon-\sigma^2\varepsilon^2}{8}V_6$ & $V_4$ & $V_5-\varepsilon(V_4+\frac{1}{2}V_6)$ & $V_6$  \\
        $V_2$ & $\frac{\sigma^2}{8}V_6+e^{\varepsilon}(V_1-\frac{\sigma^2}{8}V_6)$ & $V_2$ & $e^{-\varepsilon}V_3$ & $-\frac{1}{2}V_6+e^{\frac{\varepsilon}{2}}(V_4+\frac{1}{2}V_6)$ & $e^{-\frac{\varepsilon}{2}}V_5$ & $V_6$  \\
        $V_3$ & $V_1+\varepsilon(2V_2-\frac{1}{2}V_6)+\varepsilon^2V_3$ & $V_2+\varepsilon V_3$  & $V_3$ & $V_4+\varepsilon V_5$ & $V_5$ & $V_6$  \\
        $V_4$ & $V_1$ & $V_2-\varepsilon(\frac{1}{2}V_4+\frac{1}{4}V_6)$ & $V_3-\varepsilon V_5-\frac{\varepsilon^2}{2\sigma^2}V_6$  & $V_4$ & $V_5+\frac{\varepsilon}{\sigma^2}V_6$ & $V_6$  \\
        $V_5$ & $V_1+\varepsilon(V_4+\frac{1}{2}V_6)-\frac{\varepsilon^2}{2\sigma^2}V_6$ & $V_2+\frac{\varepsilon}{2}V_5$ & $V_3$  & $V_4-\frac{\varepsilon}{\sigma^2}V_6$ & $V_5$ & $V_6$  \\
	$V_6$ & $V_1$ & $V_2$ & $V_3$  & $V_4$ & $V_5$ & $V_6$ \\
\hline
\end{tabular}
}
\end{table}\label{tab2}

\medskip
Let $G$ represent the symmetry group of (\ref{eq1}) associated with the symmetry algebra $\mathfrak{g}$. Furthermore, let us examine an element of the six-dimensional symmetry algebra $\mathfrak{g}$, denoted by $$V=a_1V_1+a_2V_2+a_3V_3+a_4V_4+a_5V_5+a_6V_6.$$
Our objective is to simplify as many of the coefficients $a_i$ as possible through appropriate applications of adjoint maps to $V$. A significant observation here is that the function $\eta(V)=a_2^2-4a_1a_3$ is an invariant of the full adjoint action: $\eta({\rm Ad}\, g(V))=\eta (V), V\in\mathfrak{g},g\in G.$ Formulating this invariant is fundamental because it restricts how much we can simplify $V$.

\medskip
We focus on the coefficients $a_1,a_2,a_3$ of $V$ to start the classification process. If we let ${\rm Ad}(\exp(\gamma V_1))$ and ${\rm Ad}(\exp(\beta V_3))$ to act on $V$, then
\begin{equation*}
\begin{aligned}
    \tilde{V}=\sum_{i=1}^6\tilde{a}_iV_i={\rm Ad}(\exp(\gamma V_1)){\rm Ad}(\exp(\beta V_3))V
\end{aligned}
\end{equation*}
has coefficients
\begin{equation*}
\begin{aligned}
     &\tilde{a}_1=a_1-\gamma(2a_1\beta+a_2)+\gamma^2(a_1\beta^2+a_2\beta+a_3), \\
    &\tilde{a}_2=2a_1\beta+a_2-2\gamma(a_1\beta^2+a_2\beta+a_3), \\
    &\tilde{a}_3=a_1\beta^2+a_2\beta+a_3.
\end{aligned}
\end{equation*}
According to the sign of the invariant $\eta(V)$, there are three cases:

\medskip		
\textbf{\textit{Case 1.}} $\eta(V)>0$.

\medskip
We can choose $\beta=\frac{\sqrt{\eta(V)}-a_2}{2a_1}$ to be either the real root of the quadratic equation $a_1\beta^2+a_2\beta+a_3=0$, and $\gamma=a_1/(2a_1\beta+a_2)$. Then $\tilde{a}_1=\tilde{a}_3=0$ while $\tilde{a}_2=2a_1\beta+a_2=\sqrt{\eta(V)}\neq0$, so $V$ is equivalent to a multiple of $\tilde{V}=V_2+\tilde{a}_4V_4+\tilde{a}_5V_5+\tilde{a}_6V_6$ for certain scalars $\tilde{a}_4$, $\tilde{a}_5$ and $\tilde{a}_6$ depending on $a_i$. We next act on $\tilde{V}$ by ${\rm Ad}(\exp(2\tilde{a}_4V_4))$ to cancel the coefficient of $V_4$ and then by ${\rm Ad}(\exp(-2\tilde{a}_5V_5))$ to cancel the coefficient of $V_5$, leading to $\tilde{V}=V_2+\tilde{a}_6V_6$ for certain scalar $\tilde{a}_6$. Therefore, every element with $\eta(V)>0$ is equivalent to a multiple of $V_2+aV_6$ for some $a\in\mathbb{R}$.

\medskip		
\textbf{\textit{Case 2. }} $\eta(V)<0$.

\medskip
We can choose $\beta=0$, $\gamma=a_2/2a_3$, then $\tilde{a}_1=a_1-a_2^2/4a_3$, $\tilde{a}_2=0$ and $\tilde{a}_3=a_3$. Consider the action on $\tilde{V}$ by the group generated by $V_2$
\begin{eqnarray*}
    \bar{V}={\rm Ad}\left(\exp\Big(\frac{1}{2}\ln\frac{\tilde{a}_3}{\tilde{a}_1}V_2\Big)\right)\tilde{V}={\rm Ad}\left(\exp\Big(\frac{1}{2}\ln\frac{4a_3^2}{4a_1a_3-a_2^2}V_2\Big)\right)\tilde{V},
\end{eqnarray*}
which makes the coefficients of $V_1$ and $V_3$ agree, so $V$ is equivalent to a multiple of $\bar{V}=V_1+V_3+\bar{a}_4V_4+\bar{a}_5V_5+\bar{a}_6V_6$
for certain scalars $\bar{a}_4$, $\bar{a}_5$ and $\bar{a}_6$. Further use of the groups generated by $V_4$ and $V_5$ shows that
$V$ is equivalent to a multiple of $V_1+V_3+aV_6$ for some $a\in\mathbb{R}$.

\medskip		
\textbf{\textit{Case 3.}} $\eta(V)=0$.

\medskip
There are five subcases as follows.

\medskip
\textbf{\textit{Case 3.1.}} $a_1\neq0$.

\medskip
We can select $\beta=-a_2/2a_1$ to make $\tilde{a}_2=\tilde{a}_3=0$, while $\tilde{a}_1\neq0$, so $V$ is equivalent to a multiple of $\tilde{V}=V_1+\tilde{a}_4V_4+\tilde{a}_5V_5+\tilde{a}_6V_6$ for certain scalars $\tilde{a}_4$, $\tilde{a}_5$ and $\tilde{a}_6$.
Suppose $\tilde{a}_5\neq0$, then act on $\tilde{V}$ by ${\rm Ad}(\exp(\frac{2}{3}\ln|\tilde{a}_5|V_2))$ to scale the coefficients of $V_1$ and $V_5$ to $\pm1$, so $V$ is equivalent to a multiple of $\tilde{V}=V_1\pm V_5+\tilde{a}_4V_4+\tilde{a}_6V_6$. Further use of the groups generated by $V_1$ to cancel the coefficient of $V_4$ and then by $V_4$ to cancel the coefficient of $V_6$. Thus such a $V$ is equivalent to a multiple of $V_1+V_5$ or $V_1-V_5$. On the other hand, if $\tilde{a}_5=0$, then act by ${\rm Ad}(\exp(-\tilde{a}_4V_5))$ to make the coefficient of $V_4$ vanish, so $V$ is equivalent to a multiple of $V_1+aV_6$ for some $a\in\mathbb{R}$.

\medskip		
\textbf{\textit{Case 3.2.}} $a_1=a_2=0$, $a_3\neq0$.

\medskip
We can use the group generated by $V_1$ to get a nonzero coefficient of $V_1$,
which reduce to the previous case.

\medskip
\textbf{\textit{Case 3.3.}}	$a_1=a_2=a_3=0$, $a_4\neq0$.

\medskip
Without loss of generality, we can assume that $a_4=1$.
Then act on $V$ by ${\rm Ad}(\exp(-a_5V_3))$ to cancel the coefficient of $V_5$ and then by ${\rm Ad}(\exp(\sigma^2a_6V_5))$ to cancel the coefficient of $V_6$, leading to $\tilde{V}=V_4$.

\medskip
\textbf{\textit{Case 3.4.}}	$a_1=a_2=a_3=a_4=0$, $a_5\neq0$.

\medskip
We can use the group generated by $V_1$ to get a nonzero coefficient of $V_4$, which reduce to the previous case.

\medskip
\textbf{\textit{Case 3.5.}} $a_1=a_2=a_3=a_4=a_5=0$.

\medskip
We get a multiple of $V_6$, on which the adjoint representation acts trivially.

\medskip
According to the above analysis, we obtain an optimal system of one-dimensional symmetry algebras of (\ref{eq1}) as follows.
\begin{theorem}
    An optimal system of one-dimensional subalgebras of (\ref{eq1}) is provided by those generated by
\begin{eqnarray}\label{os}
    V_2+aV_6, \ \ \ V_1+V_3+aV_6, \ \ \ V_1\pm V_5, \ \ \ V_1+aV_6, \ \ \ V_4, \ \ \ V_6, \ \ \ a\in\mathbb{R}.
\end{eqnarray}
\end{theorem}
		
Once we have classified the one-dimensional subalgebras of a symmetry algebra, we can give a classification of the invariant solutions of (\ref{eq1}). The set of all invariant solutions for the one-dimensional subgroups in the optimal system (\ref{os}) will form an optimal system of invariant solutions. Moreover, any other solution invariant under a one-parameter symmetry group will correspond to a solution on this list under the group action.

\section{Symmetry reductions and exact solutions} \label{sec4}

In this section, we employ one-dimensional subalgebras of the symmetry algebra associated with (\ref{eq1}) to obtain its exact solutions through the method of symmetry reductions. For each one-parameter subgroup generated by $V_i, \ i=1,...,5$, there will be a corresponding class of invariant solutions.
		
\medskip
\textbf{\textit{Case 1.}} The invariant solutions for infinitesimal generator $V_1$ are $u=f(x)$. Substituting this expression into (\ref{eq1}),
we find the reduced equation to be
\begin{eqnarray*}
    \sigma^2f_{xx}-2kf_x+\frac{k^2}{\sigma^2}f=0.
\end{eqnarray*}
Hence the invariant solutions are
\begin{eqnarray}\label{inv1}
    u(t,x)=f(x)=(c_1+c_2x)e^{\frac{k}{\sigma^2}x},
\end{eqnarray}
in which $c_1$ and $c_2$ are the constants of integration.

\medskip
\textbf{\textit{Case 2.}} The invariant solutions corresponding to infinitesimal generator $V_2$ are $$u=x^{\frac{1}{2}}e^{\frac{k}{\sigma^2}x-\frac{\sigma^2}{8}t} f(y),$$  where $y=t^{-\frac{1}{2}}\ln x$. Performing the symmetry reduction, we arrive at the following reduced equation
\begin{eqnarray*}
    \sigma^2f_{yy}+yf_y=0.
\end{eqnarray*}
Solving the above equation yields the invariant solutions
\begin{eqnarray}\label{inv2}
    u(t,x)=x^{\frac{1}{2}}e^{\frac{k}{\sigma^2}x-\frac{\sigma^2}{8}t}\left(c_1\int e^{-\frac{y^2}{2\sigma^2}}{\rm d}y\Big{|}_{y=t^{-\frac{1}{2}}\ln x}+c_2\right),
\end{eqnarray}
where $c_1$ and $c_2$ are arbitrary constants. We set $k=0.5$, $c_1=2$ and $c_2=-1$. Figure 1 displays the dynamic behavior of the pricing formula (\ref{inv2}) under different $\sigma$ values. When the volatility $\sigma$ is small, the price is more sensitive to the state
$x$ (initial value of the GMR process). Specifically, at a fixed time $t$, as the initial value $x$ increases from 1 to 5, the price changes more significantly for smaller $\sigma$. Additionally, prices are higher for the same $x$ when $\sigma$ is smaller. Conversely, when
the volatility $\sigma$ is large, the price is more sensitive to time $t$. Over time, the price tends to stabilize.
\begin{figure}[H]
	\centering
	\subfigure[$\sigma = 1$]{
		\begin{minipage}[b]{0.3\textwidth}
			\includegraphics[width=1\textwidth]{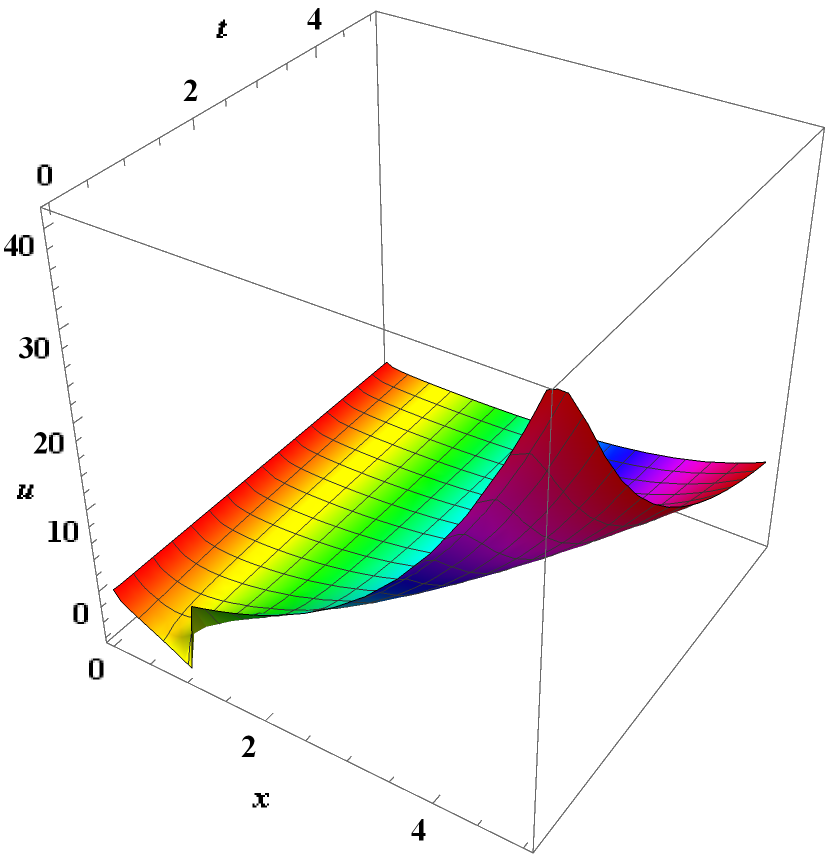}
		\end{minipage}
		\label{figsub1}
	}
    \subfigure[$\sigma = 2$]{
    	\begin{minipage}[b]{0.3\textwidth}
   		\includegraphics[width=1\textwidth]{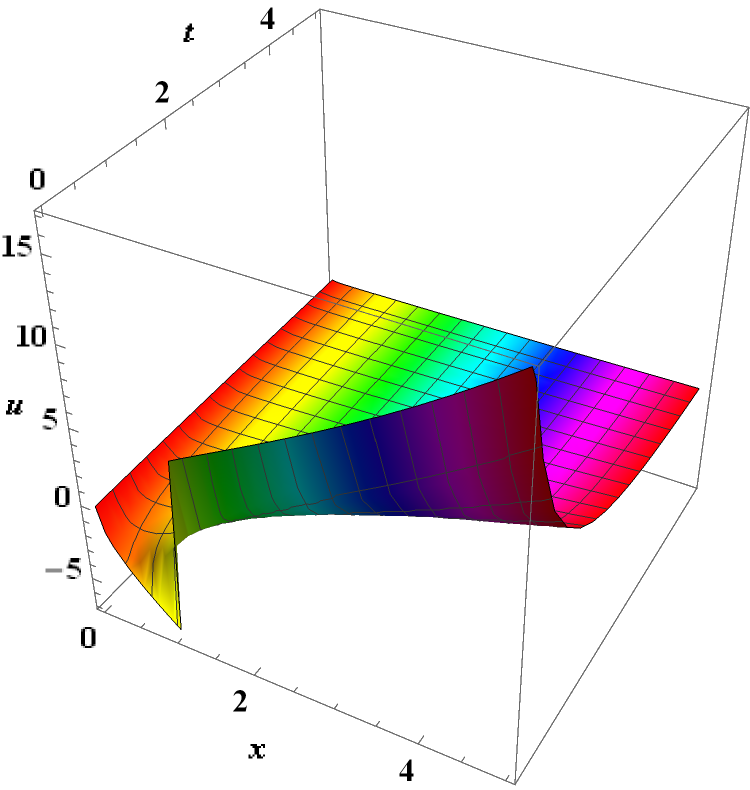}
    	\end{minipage}
	\label{figsub2}
    }
    \subfigure[$\sigma = 3$]{
    	\begin{minipage}[b]{0.3\textwidth}
   		\includegraphics[width=1\textwidth]{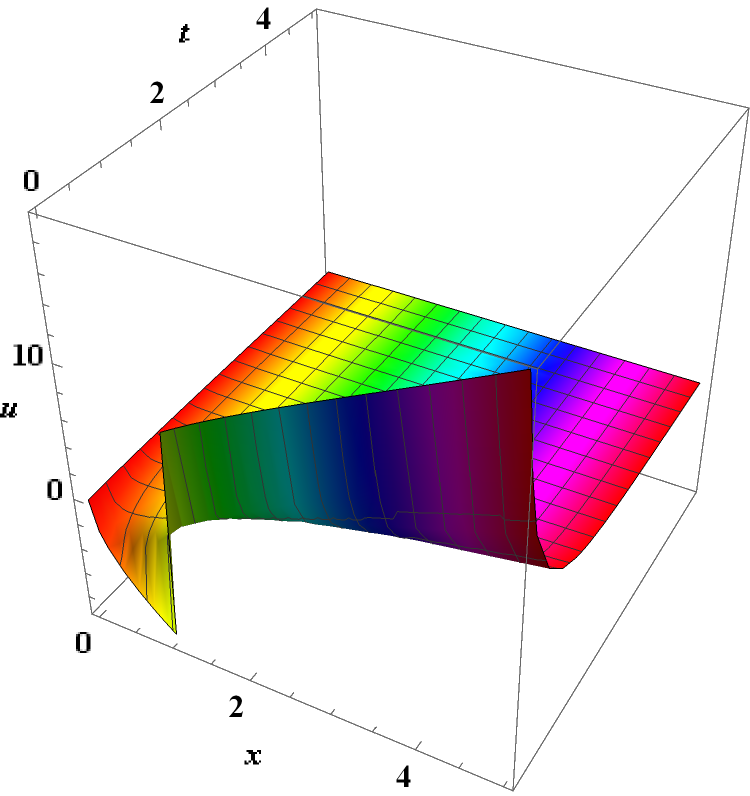}
    	\end{minipage}
	\label{figsub3}
    }
	\caption{Three different dynamical structures of (\ref{inv2}) for $k = 0.5$, $c_1 = 2$, $c_2 = -1$.}
	\label{fig1}
\end{figure}
		
\textbf{\textit{Case 3.}} The invariant solutions corresponding to $V_3$ are given by
\begin{equation*}
    u=x^{\frac{1}{2}}(\ln x)^{-\frac{1}{2}}e^{\frac{k}{\sigma^2}x-\frac{\sigma^2}{8}t-\frac{1}{2\sigma^2}t^{-1}\ln^2x}f(z),
\end{equation*}
where $z=t^{-1}\ln x$. The corresponding reduced equation is
\begin{eqnarray*}\label{inv12}
    z^2f_{zz}-zf_z+\frac{3}{4}f=0,
\end{eqnarray*}
and the invariant solutions are
\begin{eqnarray}\label{inv3}
    u(t,x)=t^{-\frac{1}{2}}x^{\frac{1}{2}}e^{\frac{k}{\sigma^2}x-\frac{\sigma^2}{8}t-\frac{1}{2\sigma^2}t^{-1}\ln^2x}(c_1+c_2t^{-1}\ln x),
\end{eqnarray}
where $c_1$ and $c_2$ are arbitrary constants. Set $k=1$, $c_1=2$ and $c_2=-1$. Figure 2 shows the dynamic behavior of the pricing formula (\ref{inv3}) under different $\sigma$ values. As the volatility $\sigma$ increases, the range of $x$ where the price oscillates intensely also increases, with denser oscillations and greater amplitude as $x$ increases. For smaller $\sigma$ and sufficiently large $x$, the price initially decreases to a minimum point and then rapidly increases to a maximum point before gradually stabilizing.
\begin{figure}[H]
	\centering
	\subfigure[$\sigma = 1$]{
		\begin{minipage}[b]{0.3\textwidth}
			\includegraphics[width=1\textwidth]{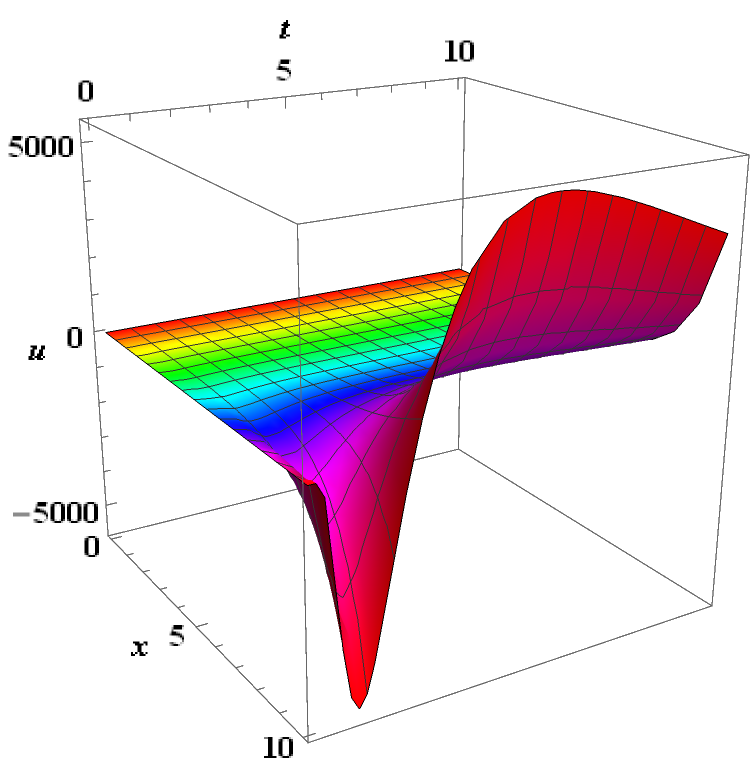}
		\end{minipage}
		\label{figsub21}
	}
    \subfigure[$\sigma = 2$]{
    	\begin{minipage}[b]{0.3\textwidth}
   		\includegraphics[width=1\textwidth]{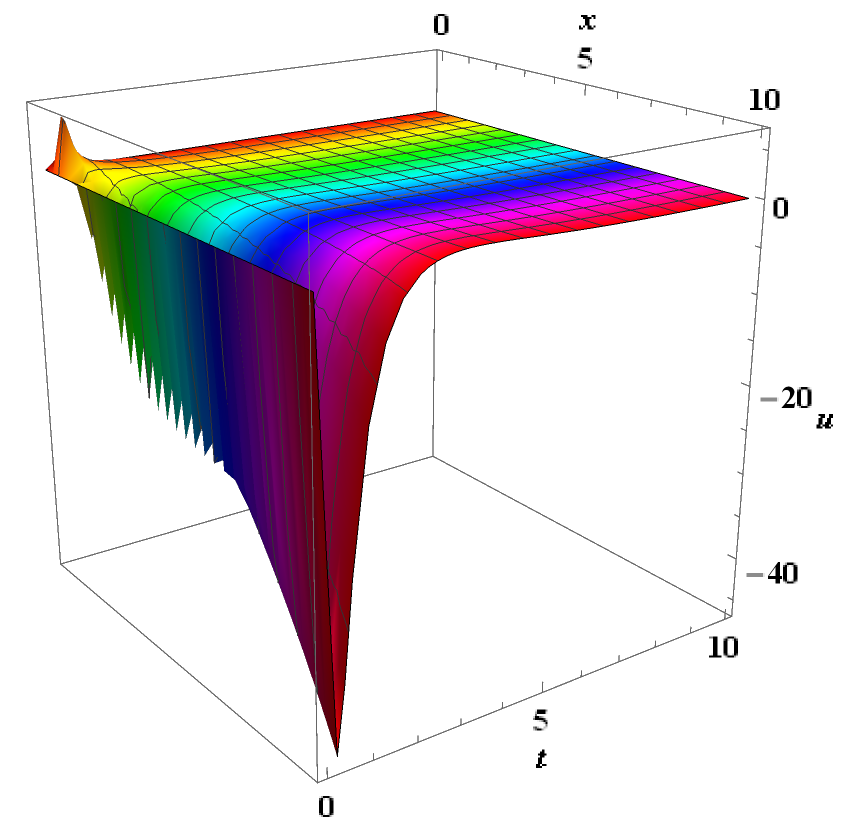}
    	\end{minipage}
	\label{figsub22}
    }
    \subfigure[$\sigma = 3$]{
    	\begin{minipage}[b]{0.3\textwidth}
   		\includegraphics[width=1\textwidth]{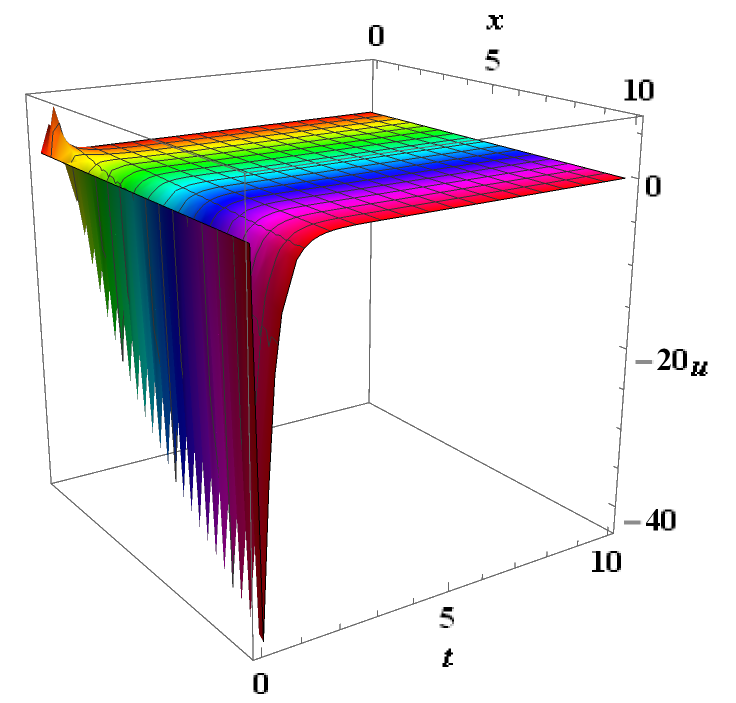}
    	\end{minipage}
	\label{figsub23}
    }
	\caption{Three different dynamical structures of (\ref{inv3}) for $k = 1$, $c_1 = 2$, $c_2 = -1$.}
	\label{fig2}
\end{figure}

\medskip		
\textbf{\textit{Case 4.}} The invariant solutions for $V_4$ are $u=f(t)e^{\frac{k}{\sigma^2}x}$. By substituting it into (\ref{eq1}), we obtain the invariant solutions
\begin{eqnarray}\label{c4}
    u=ce^{\frac{k}{\sigma^2}x},
\end{eqnarray}
where $c$ is an arbitrary constant, which we earlier found as the invariant solutions corresponding to $V_1$. This means that a given solution may be invariant under more than one subgroup of the full symmetry group.
	
\medskip	
\textbf{\textit{Case 5.}} For infinitesimal generator $V_5$, we derive the invariant solutions of the following type $$u=f(t)x^{\frac{1}{2}}e^{\frac{k}{\sigma^2}x -\frac{1}{2\sigma^2}t^{-1}\ln^2x},$$
where the invariant function $f(t)$ satisfies the reduced equation
\begin{eqnarray*}
    f_t+\left(\frac{\sigma^2}{8}+\frac{1}{2t}\right)f=0.
\end{eqnarray*}
Solving the reduced equation, we derive the invariant solutions
\begin{eqnarray}\label{inv4}
    u(t,x)=ct^{-\frac{1}{2}}x^{\frac{1}{2}}e^{\frac{k}{\sigma^2}x-\frac{\sigma^2}{8}t-\frac{1}{2\sigma^2}t^{-1}\ln^2x},
\end{eqnarray}
where $c$ is an arbitrary constant. Figure 3 illustrates the dynamic behavior of the pricing formula (\ref{inv4}) under different $\sigma$ values. For low volatility and large $x$, the price initially increases to a maximum point and then slowly decreases, stabilizing over time. For high volatility, the price exhibits more intense fluctuations with a broader range of variation.
\begin{figure}[H]
	\centering
	\subfigure[$\sigma = 1$]{
		\begin{minipage}[b]{0.3\textwidth}
			\includegraphics[width=1\textwidth]{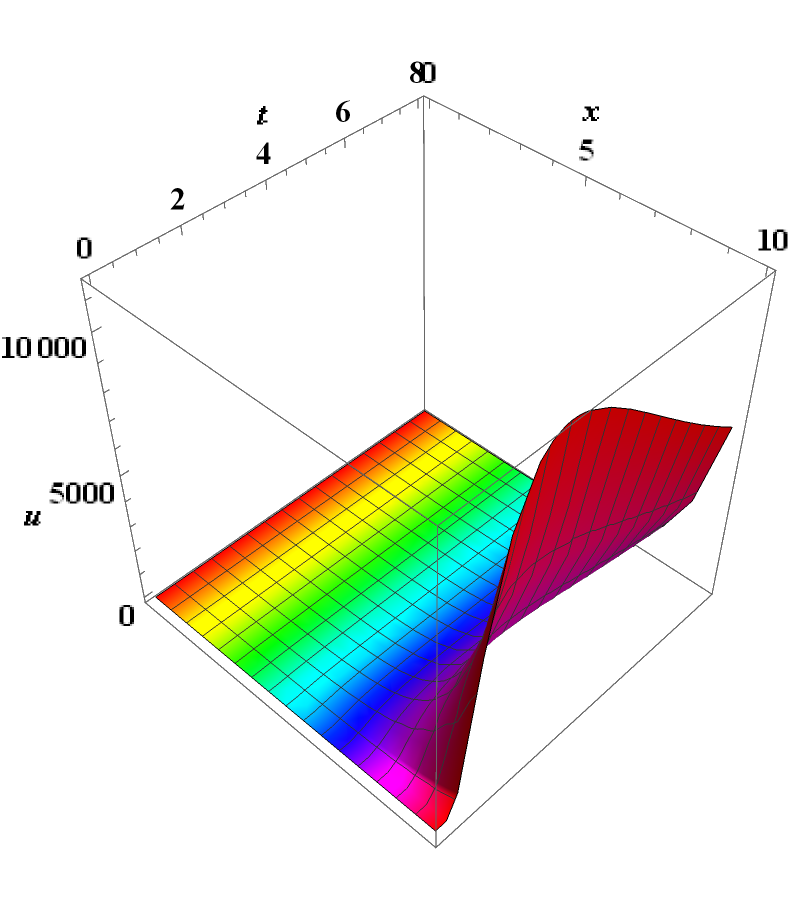}
		\end{minipage}
		\label{figsub41}
	}
    \subfigure[$\sigma = 2$]{
    	\begin{minipage}[b]{0.3\textwidth}
   		\includegraphics[width=1\textwidth]{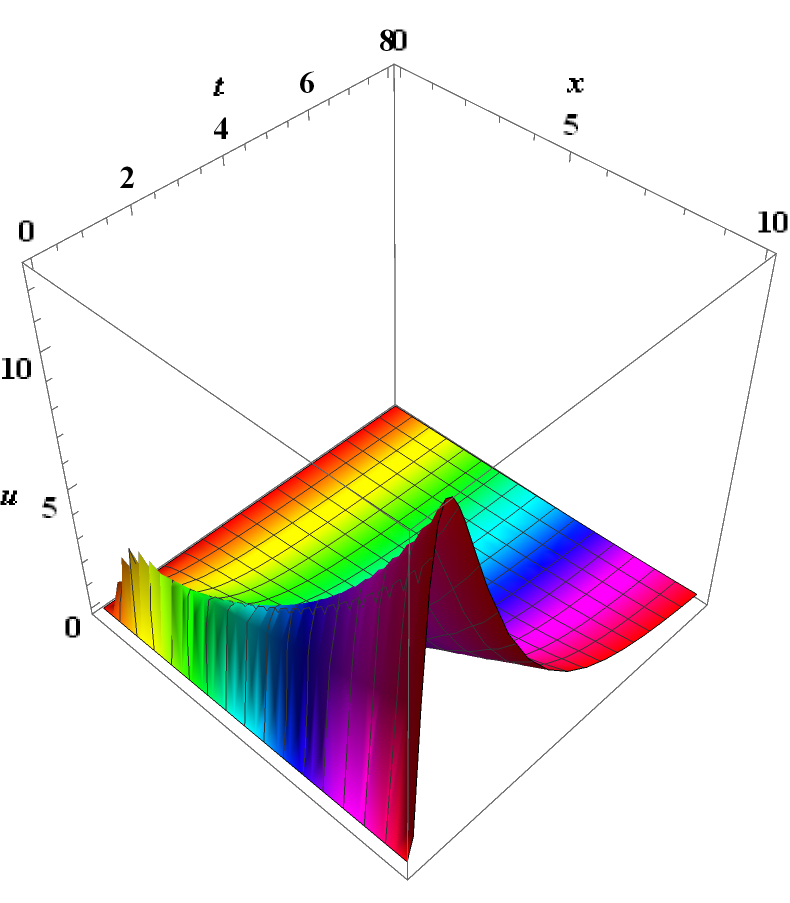}
    	\end{minipage}
	\label{figsub42}
    }
    \subfigure[$\sigma = 3$]{
    	\begin{minipage}[b]{0.3\textwidth}
   		\includegraphics[width=1\textwidth]{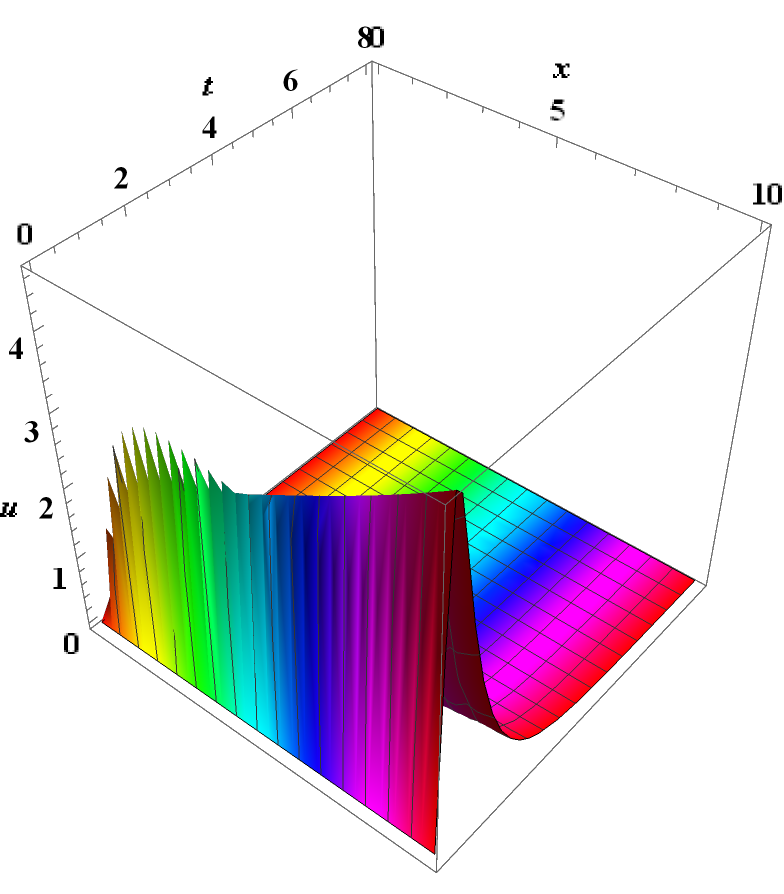}
    	\end{minipage}
	\label{figsub43}
    }
	\caption{Three different dynamical structures of (\ref{inv4}) for $k = 1$, $c = 1$.}
	\label{fig4}
\end{figure}

We can clearly expand this list of invariant solutions by considering additional one-parameter subgroups obtained from more general linear combinations of the infinitesimal generators of the full symmetry group. However, without some means of classifying these solutions, it is somewhat pointless to continue. Once we establish an appropriate classification procedure, we will be able to identify the most comprehensive set of invariant solutions.
	
\medskip	
Indeed, according to our optimal system (\ref{os}), we only need to find the invariant solutions for one-parameter subgroups generated by (a) $V_2+aV_6$; (b) $V_1+V_3+aV_6$; (c) $V_1\pm V_5$; (d) $V_1+aV_6$.
		
(a) The invariant solutions corresponding to $V_2+aV_6$ are given by
\begin{eqnarray*}
    u=t^ax^{\frac{1}{2}}e^{\frac{k}{\sigma^2}x-\frac{\sigma^2}{8}t}f(y),
\end{eqnarray*}
where $y=t^{-\frac{1}{2}}\ln x$. The substitution of this expression into (\ref{eq1}) leads to the reduced equation
\begin{eqnarray*}
    \sigma^2f_{yy}+yf_y-2af=0.
\end{eqnarray*}
The solutions of this equation can be written in terms of parabolic cylinder functions
\begin{eqnarray*}
    f(y)=e^{-\frac{1}{4\sigma^2}y^2}\left(c_1U\Big(2a+\frac{1}{2},\frac{y}{\sigma}\Big)+c_2V\Big(2a+\frac{1}{2},\frac{y}{\sigma}\Big)\right),
\end{eqnarray*}
where $U(\cdot,\cdot)$ and $V(\cdot,\cdot)$ are parabolic cylinder functions \cite{abramowitz1968handbook}, $c_1$ and $c_2$ are constants. Thus the invariant solutions are
\begin{equation}\label{inv5}
    \begin{aligned}
    u(t,x)=
    &t^ax^{\frac{1}{2}}e^{\frac{k}{\sigma^2}x-\frac{\sigma^2}{8}t-\frac{1}{4\sigma^2}t^{-1}\ln^2x}\\
    &\cdot\left(c_1U\Big(2a+\frac{1}{2},
    \frac{1}{\sigma}t^{-\frac{1}{2}}\ln x\Big)+c_2V\Big(2a+\frac{1}{2},\frac{1}{\sigma}t^{-\frac{1}{2}}\ln x\Big)\right).
    \end{aligned}
\end{equation}
Set $a=2$, $k=1$, $c_1=2$ and $c_2=1$. Figure 4 shows the dynamic behavior of the pricing formula (\ref{inv5}) under different $\sigma$ values. When the volatility $\sigma$ is large, the price curve near the state $x=1$ becomes steeper over time, while the price curve farther from $x=1$ shows more gradual changes. At the same time, the peak value of the price
is lower.
\begin{figure}[H]
\centering
\subfigure[$\sigma = 1$]{
		\begin{minipage}[b]{0.3\textwidth}
			\includegraphics[width=1\textwidth]{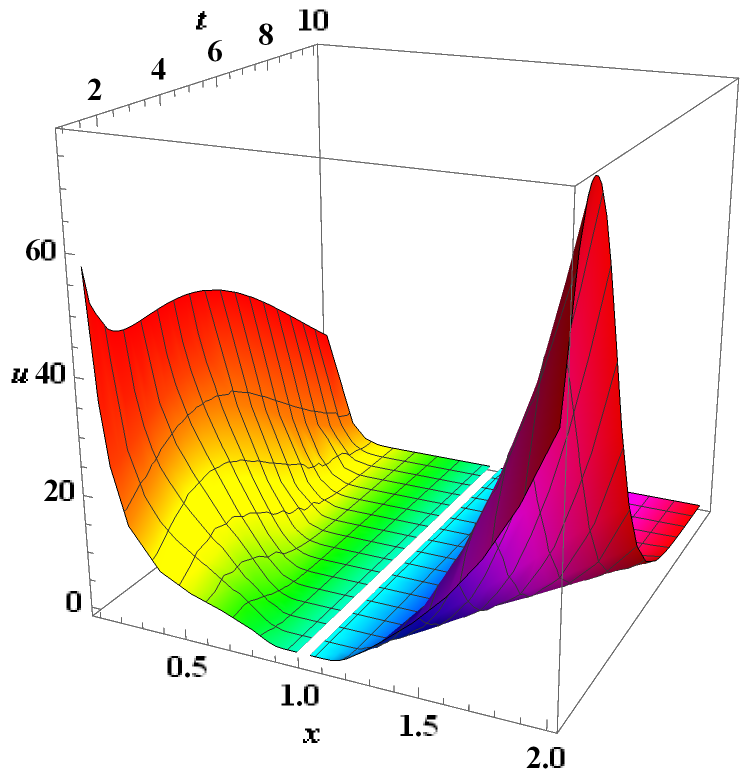}
		\end{minipage}
		\label{figsuba1}
	}
    \subfigure[$\sigma = 2$]{
    	\begin{minipage}[b]{0.307\textwidth}
   		\includegraphics[width=1\textwidth]{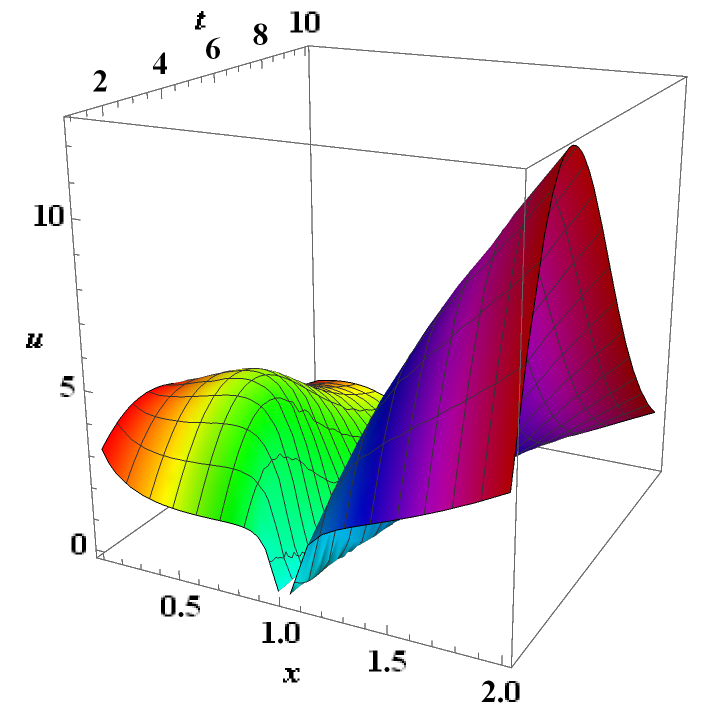}
    	\end{minipage}
	\label{figsuba2}
    }
    \subfigure[$\sigma = 3$]{
    	\begin{minipage}[b]{0.3125\textwidth}
   		\includegraphics[width=1\textwidth]{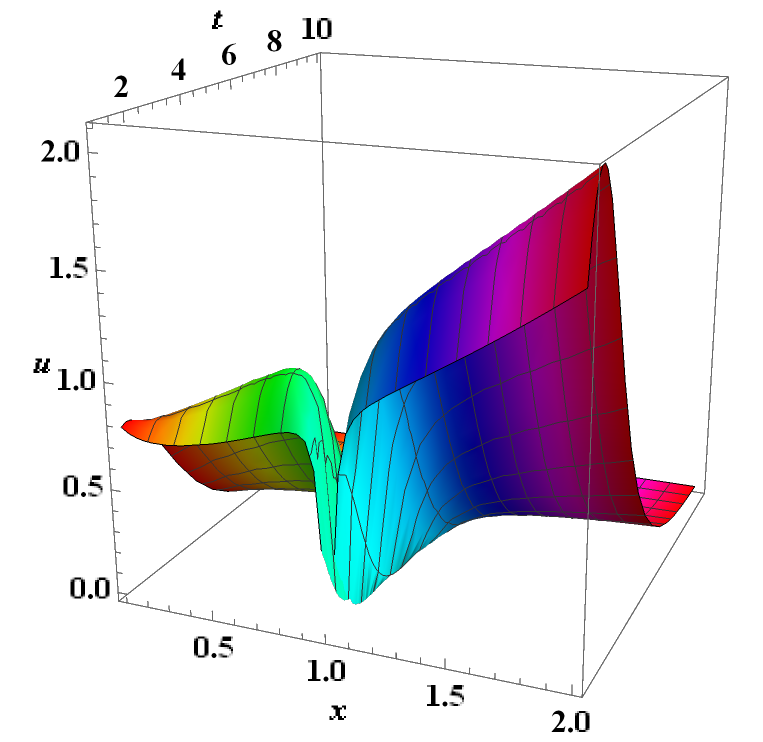}
    	\end{minipage}
	\label{figsuba3}
    }
\caption{Three different dynamical structures of (\ref{inv5}) for $a = 2$, $k = 1$, $c_1 = 2$, $c_2 = 1$.}
	\label{figa}
\end{figure}
		
(b) Consider the linear combination $V_1+V_3+aV_6$, the corresponding invariant solutions are
\begin{eqnarray*}
    u=(1+t^2)^{-\frac{1}{4}}x^{\frac{1}{2}}e^{\frac{k}{\sigma^2}x-\frac{\sigma^2}{8}t+\frac{\sigma^2+8a}{8}\arctan t-\frac{1}{2\sigma^2}t(1+t^2)^{-1}\ln^2x}f(s),
\end{eqnarray*}
where $s=(1+t^2)^{-\frac{1}{2}}\ln x$, and the reduced equation is
\begin{eqnarray*}
    \sigma^2f_{ss}+\left(\frac{1}{\sigma^2}s^2-\frac{\sigma^2+8a}{4}\right)f=0.
\end{eqnarray*}
The solutions of this linear ordinary differential equation can be written as
\begin{eqnarray*}
    f(s)=c_1W\Big(\frac{\sigma^2+8a}{8},\frac{s}{\sqrt{2}\sigma}\Big)+c_2W\Big(\frac{\sigma^2+8a}{8},-\frac{s}{\sqrt{2}\sigma}\Big),
\end{eqnarray*}
where $W(\cdot,\cdot)$ is the parabolic cylinder function \cite{abramowitz1968handbook} with imaginary arguments,  $c_1$ and $c_2$ are arbitrary constants. Hence the invariant solutions are
\begin{equation}\label{inv6}
    \begin{aligned}
    u(t,x)
    =~&(1+t^2)^{-\frac{1}{4}}x^{\frac{1}{2}}e^{\frac{k}{\sigma^2}x-\frac{\sigma^2}{8}t
    +\frac{\sigma^2+8a}{8}\arctan t-\frac{1}{2\sigma^2}t(1+t^2)^{-1}\ln^2x} \\
    &\! \cdot \! \left(\! c_1W \! \Big(\! \frac{\sigma^2 \! + \! 8a}{8},\frac{(1 \!+\! t^2)^{-\frac{1}{2}}\ln x}{\sqrt{2}\sigma} \! \Big)
    \!+\! c_2 W \! \Big(\! \frac{\sigma^2 \! +\! 8a}{8},-\frac{(1\! +\! t^2)^{-\frac{1}{2}}\ln x}{\sqrt{2}\sigma} \! \Big)\!\right)\!.\!
\end{aligned}
\end{equation}
		
(c) The invariant solutions corresponding to $V_1+V_5$ are
\begin{equation*}
   u=e^{\frac{k}{\sigma^2}x+\frac{1}{4}t^2+\frac{1}{3\sigma^2}t^3-\frac{1}{\sigma^2}t\ln x}f(v),
\end{equation*}
where $v=\ln x-\frac{1}{2}t^2$. Substituting this expression into (\ref{eq1}), we obtain the reduced equation
\begin{eqnarray*}
    \sigma^2f_{vv}-\sigma^2f_v+\frac{2}{\sigma^2}vf=0.
\end{eqnarray*}
The general solutions of this equation may be expressed in terms of Airy functions
\begin{eqnarray*}
    f(v)=e^{\frac{1}{2}v}\left(c_1\textrm{Ai}\Big(\big(4\sigma\big)^{-\frac{4}{3}}\big(\sigma^4-8v\big)\Big)
    +c_2\textrm{Bi}\Big(\big(4\sigma\big)^{-\frac{4}{3}}\big(\sigma^4-8v\big)\Big)\right),
\end{eqnarray*}
where $\textrm{Ai}(\cdot)$ and $\textrm{Bi}(\cdot)$ are Airy functions \cite{abramowitz1968handbook}, $c_1$ and $c_2$ are arbitrary constants. Thus the invariant solutions are
\begin{equation}\label{inv7-1}
\begin{aligned}
     u(t,x)
     =~&x^{\frac{1}{2}}e^{\frac{k}{\sigma^2}x+\frac{1}{3\sigma^2}t^3-\frac{1}{\sigma^2}t\ln x}\\
     &\! \cdot \! \left(\! c_1\textrm{Ai}\Big(\! \big(4\sigma\big)^{-\frac{4}{3}}\big(\sigma^4 \! - \! 8\ln x \! + \! 4t^2\big)\Big)
     \! +\! c_2\textrm{Bi}\Big(\! \big(4\sigma\big)^{-\frac{4}{3}}\big(\sigma^4 \! - \! 8\ln x \! + \! 4t^2\big)\! \Big)\!\right)\!.\!
\end{aligned}
\end{equation}
Consider $k=1$, $c_1=-1$ and $c_2=1$. Figure 5 demonstrates the dynamic behavior of the pricing formula (\ref{inv7-1}) under different $\sigma$ values. As the volatility $\sigma$ increases, the price curve experiences rapid increases after a period. The larger the $\sigma$, the faster the price rises.
\begin{figure}[H]
\centering
\subfigure[$\sigma = 1$]{
	\begin{minipage}[b]{0.3\textwidth}
		\includegraphics[width=1\textwidth]{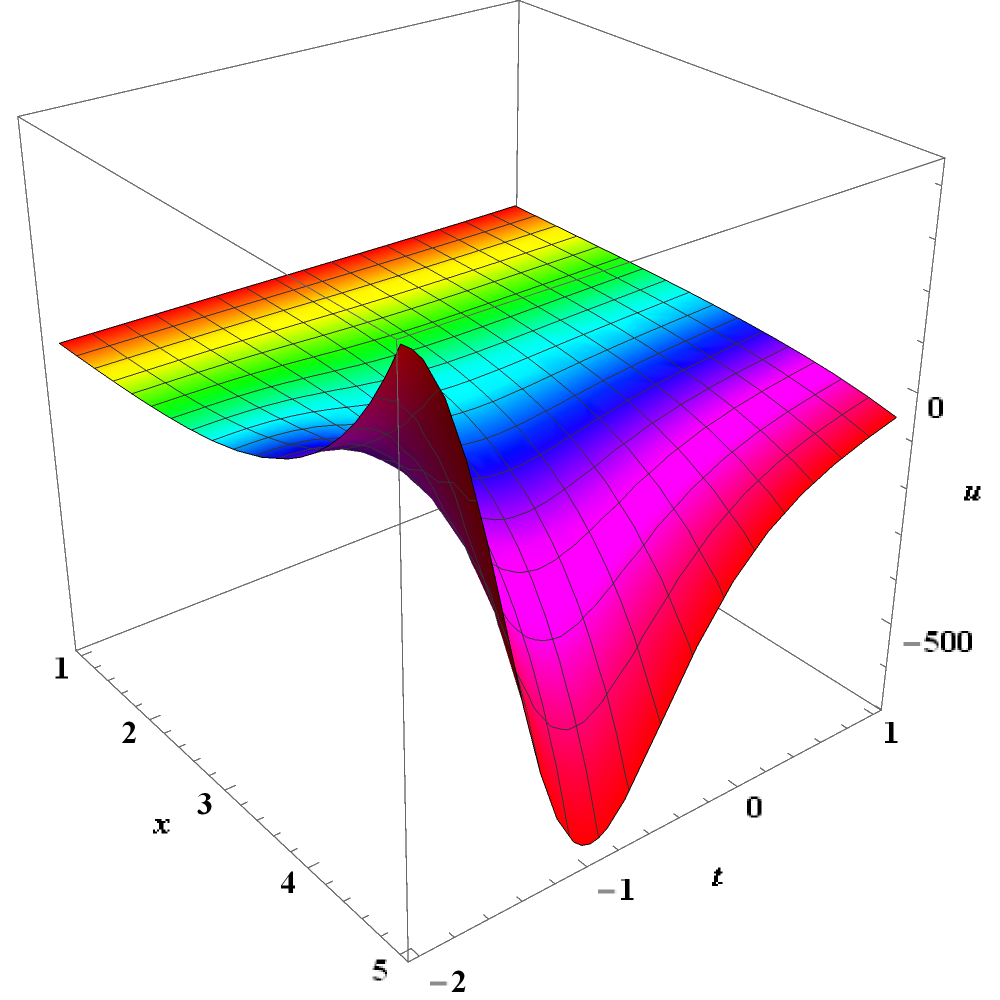}
	\end{minipage}
	\label{figsubc11}
}
\subfigure[$\sigma = 2$]{
	\begin{minipage}[b]{0.3\textwidth}
		\includegraphics[width=1\textwidth]{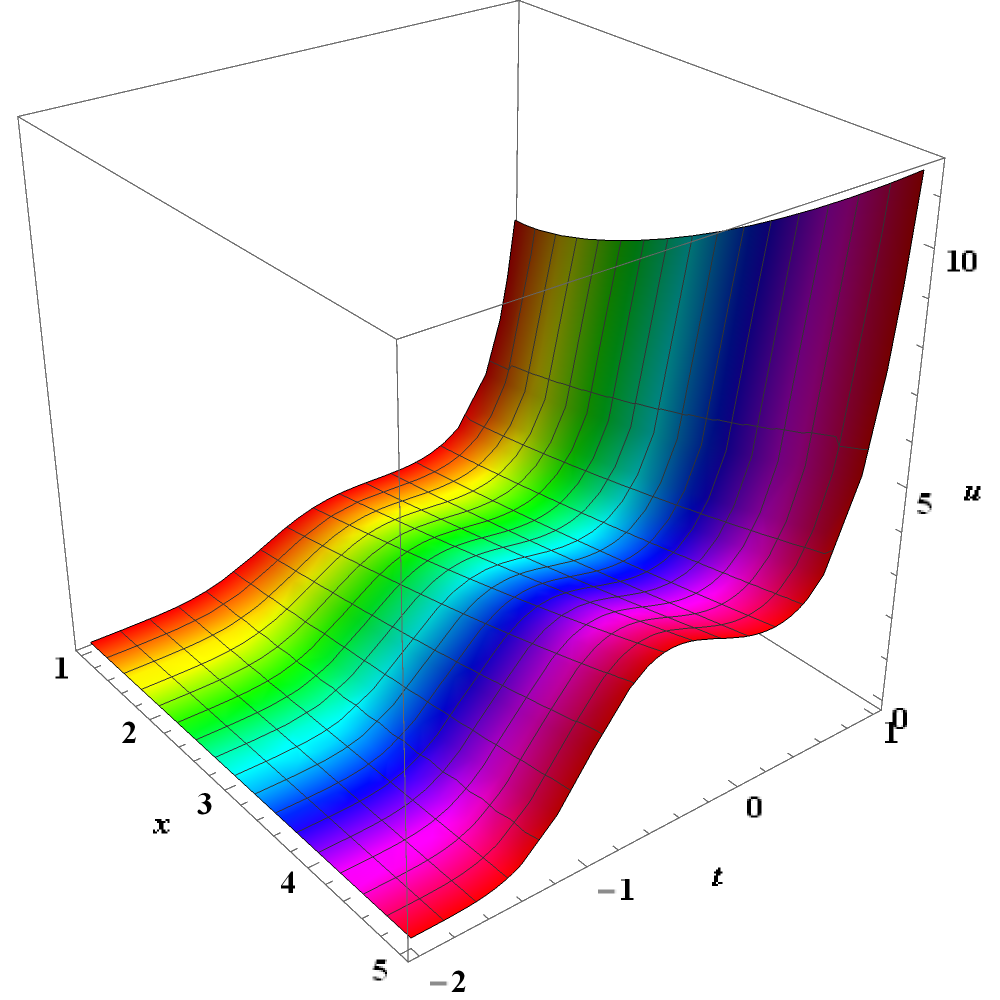}
	\end{minipage}
	\label{figsubc12}
}
\subfigure[$\sigma = 3$]{
	\begin{minipage}[b]{0.3\textwidth}
		\includegraphics[width=1\textwidth]{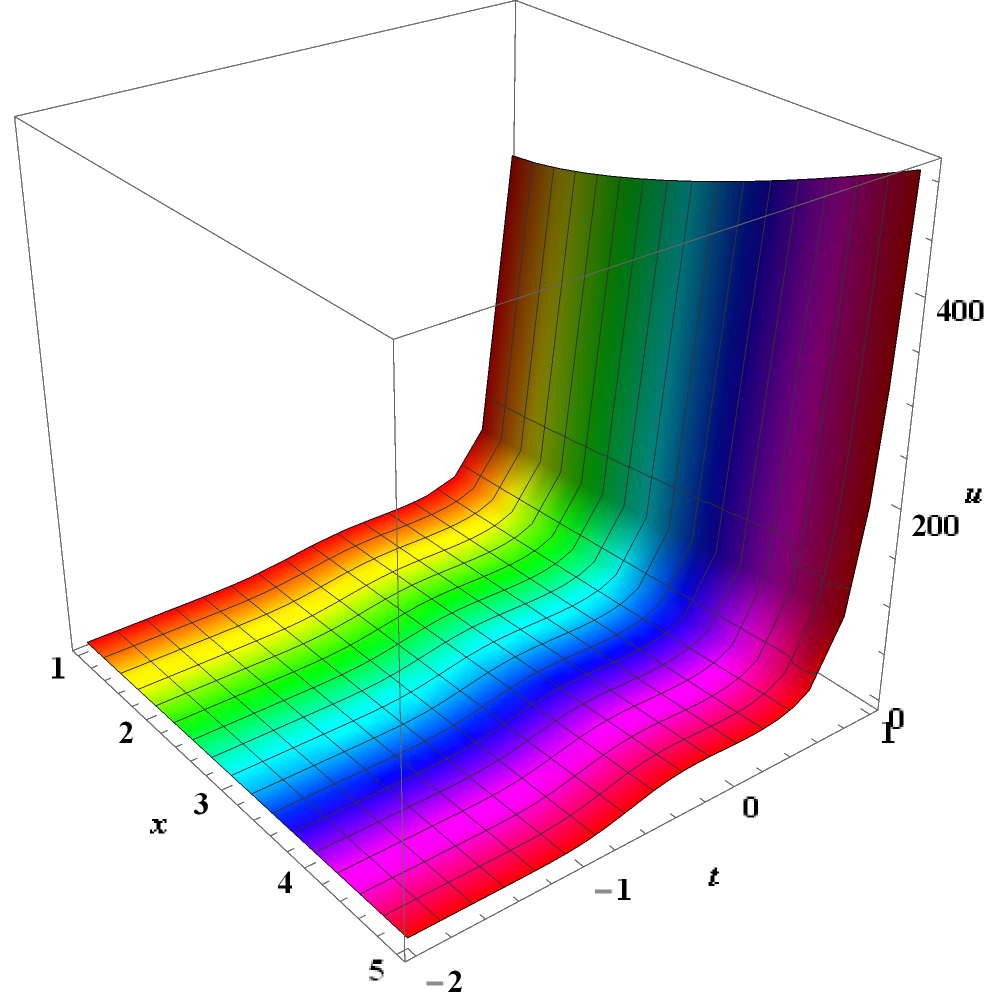}
	\end{minipage}
	\label{figsubc13}
}
\caption{Three different dynamical structures of (\ref{inv7-1}) for $c_1 =-1$, $c_2 = 1$, $k = 1$.}
\label{figc1}
\end{figure}

Similarly, the invariant solutions corresponding to $V_1-V_5$ are given by
\begin{equation}\label{inv7-2}
\begin{aligned}
    u(t,x)=
    &x^{\frac{1}{2}}e^{\frac{k}{\sigma^2}x+\frac{1}{3\sigma^2}t^3+\frac{1}{\sigma^2}t\ln x}\\
    &\!\cdot \!\left(\! c_1\textrm{Ai}\Big(\!\big(4\sigma\big)^{-\frac{4}{3}}\big(\sigma^4 \! + \! 8\ln x \! + \! 4t^2\big)\!\Big)
    +c_2\textrm{Bi}\Big(\!\big(4\sigma\big)^{-\frac{4}{3}}\big(\sigma^4 \! + \! 8\ln x \! + \! 4t^2\big) \! \Big)\!\right)\!.\!
\end{aligned}
\end{equation}
For $k=1$, $c_1=2$ and $c_2=1$. Figure 6 illustrates the dynamic behavior of the pricing formula (\ref{inv7-2}) under different $\sigma$ values. When the volatility $\sigma$ is small, the price is more sensitive to changes in time $t$ and state $x$. As $\sigma$ increases, the growth rate of the price first decreases and then increases. At an optimal $\sigma$, the price has the minimum growth rate.
\begin{figure}[H]
\centering
\subfigure[$\sigma = 1$]{
		\begin{minipage}[b]{0.3\textwidth}
			\includegraphics[width=1\textwidth]{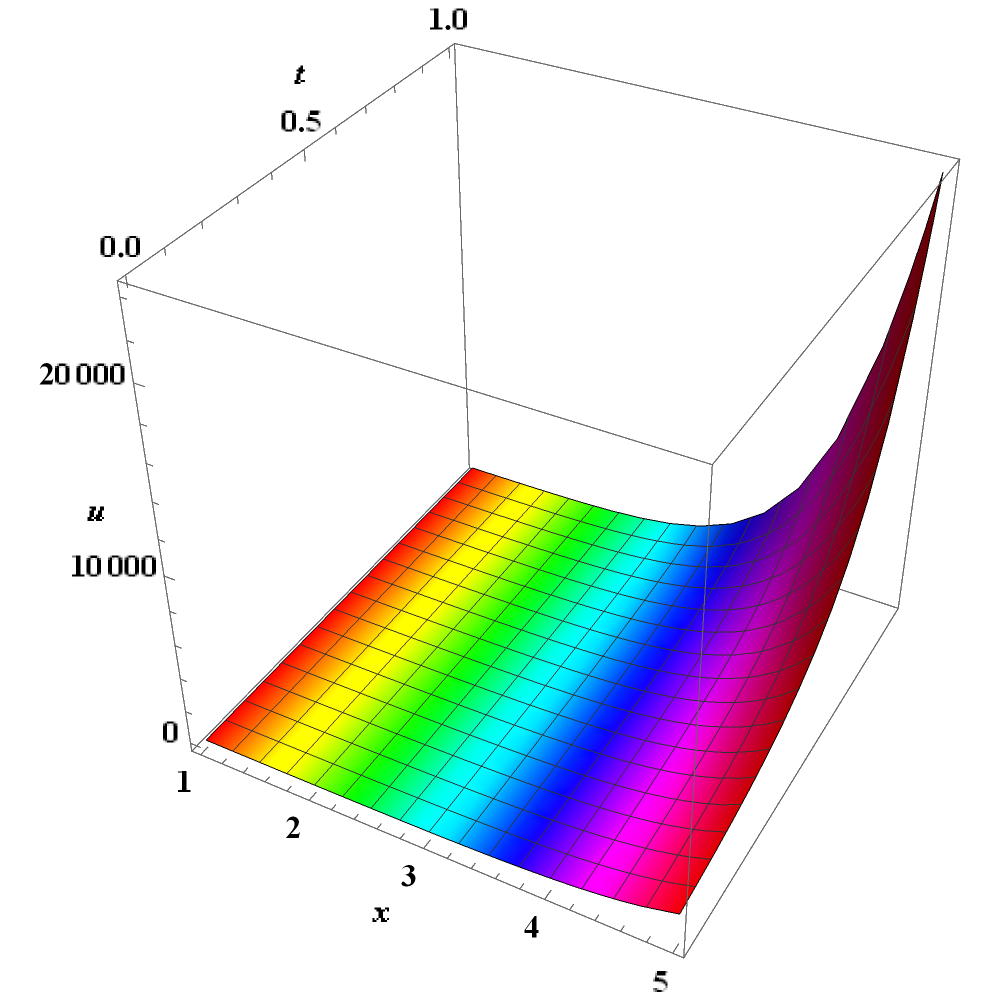}
		\end{minipage}
		\label{figsubc1}
	}
    \subfigure[$\sigma = 2$]{
    	\begin{minipage}[b]{0.3\textwidth}
   		\includegraphics[width=1\textwidth]{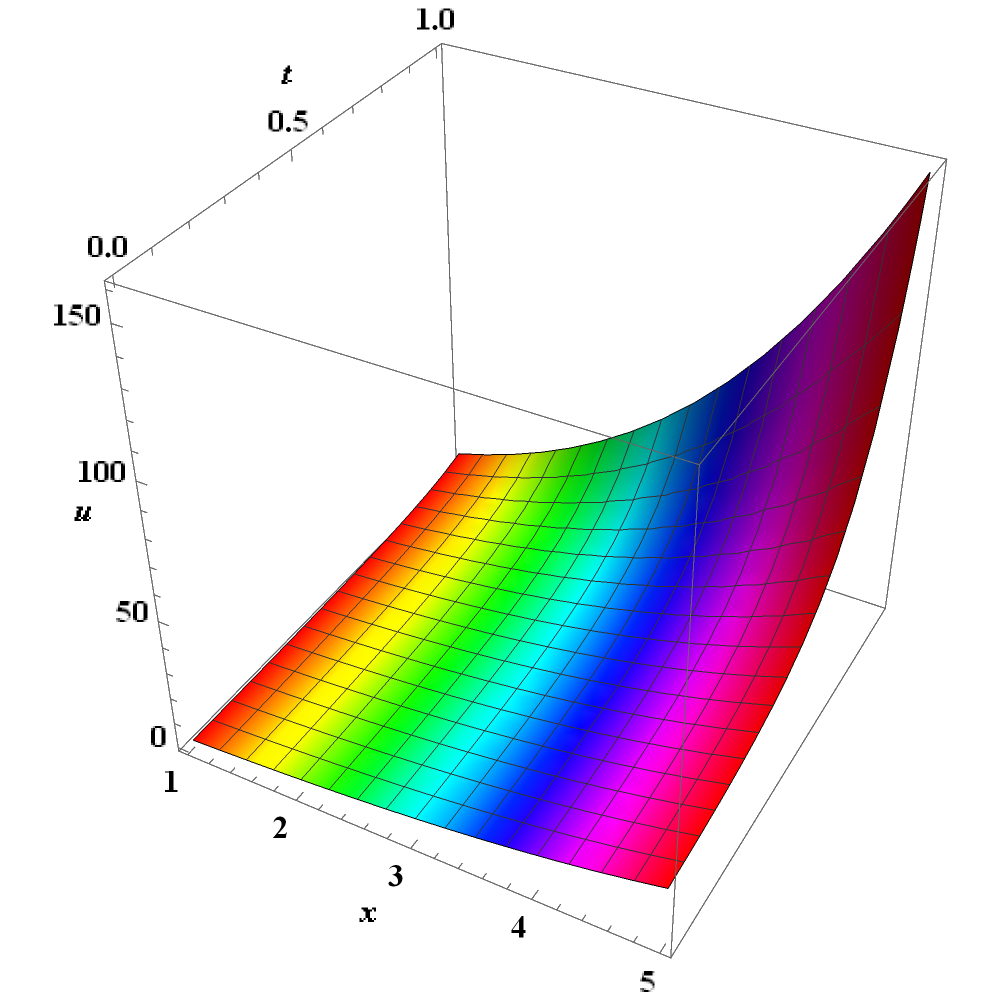}
    	\end{minipage}
	\label{figsubc2}
    }
    \subfigure[$\sigma = 3$]{
    	\begin{minipage}[b]{0.3\textwidth}
   		\includegraphics[width=1\textwidth]{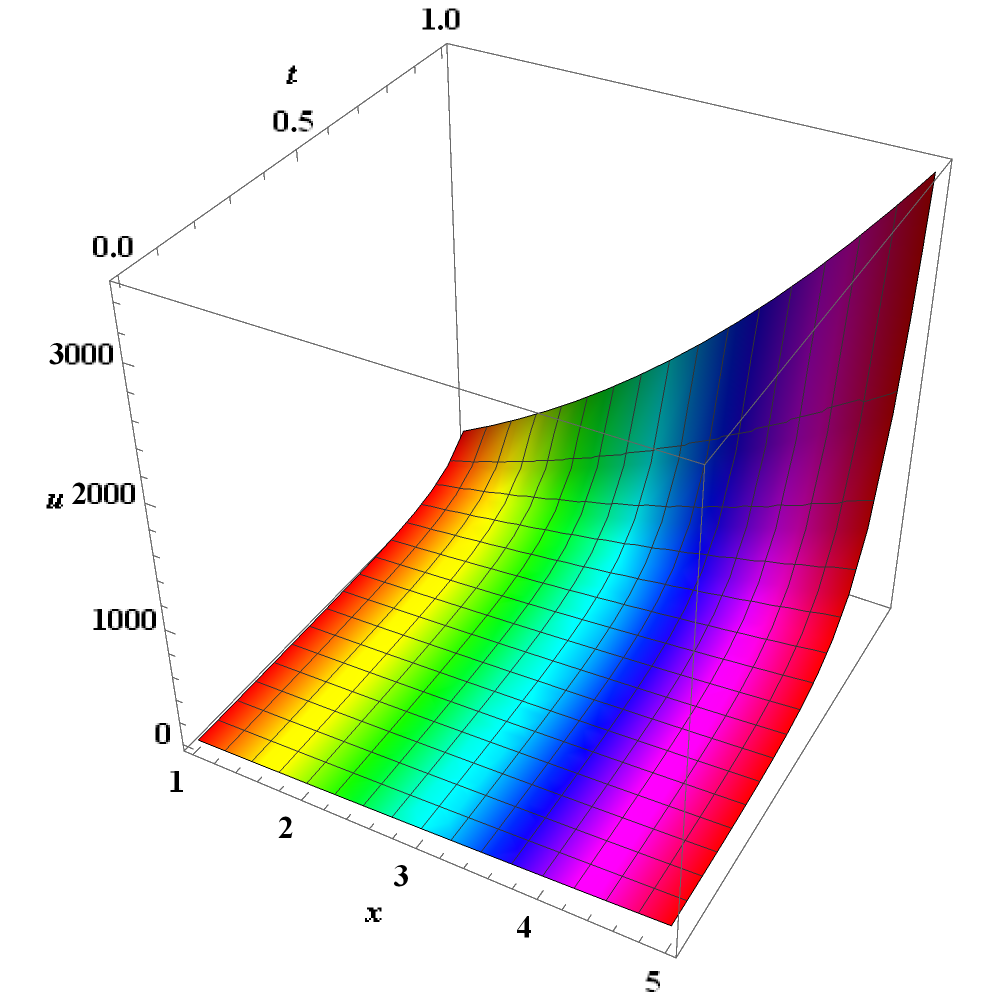}
    	\end{minipage}
	\label{figsubc3}
    }
    \caption{Three different dynamical structures of (\ref{inv7-2}) for $c_1 =2$, $c_2 = 1$, $k = 1$.}
    \label{figc}
\end{figure}
		
(d) The invariant solutions for $V_1+aV_6$ are $u=f(x)e^{at}$. Upon substituting this expression into (\ref{eq1}), we immediately obtain the reduced equation
\begin{eqnarray*}
    \sigma^2x^2f_{xx}-2kx^2f_x+\frac{k^2}{\sigma^2}x^2f=2af,
\end{eqnarray*}
which can be solved to yield
\begin{eqnarray*}
    f(x)=\left\{\begin{array}{lll}x^{\frac{1}{2}}e^{\frac{k}{\sigma^2}x}\left(c_1x^{\frac{1}{2\sigma}\sqrt{\sigma^2+8a}}
    +c_2x^{-\frac{1}{2\sigma}\sqrt{\sigma^2+8a}}\right),\qquad \qquad \qquad \ &\!\!\!\!\sigma^2+8a>0, \\
    x^{\frac{1}{2}}e^{\frac{k}{\sigma^2}x}(c_1+c_2\ln x),\qquad \qquad \qquad \qquad \qquad \qquad \qquad \qquad &\!\!\!\!\sigma^2+8a=0, \\
    x^{\frac{1}{2}}e^{\frac{k}{\sigma^2}x}\left(c_1\cos(\frac{\sqrt{-\sigma^2-8a}}{2\sigma}\ln x)
    +c_2\sin(\frac{\sqrt{-\sigma^2-8a}}{2\sigma}\ln x)\right),\ &\!\!\!\!\sigma^2+8a<0, \end{array}\right.
\end{eqnarray*}
where $c_1$ and $c_2$ are arbitrary constants. Hence the invariant solutions are
\begin{eqnarray}\label{inv8}
			u \!= \!\left\{\! \begin{array}{lll}x^{\frac{1}{2}}e^{\frac{k}{\sigma^2}x+at}\left(c_1x^{\frac{1}{2\sigma}\sqrt{\sigma^2+8a}}
				+c_2x^{-\frac{1}{2\sigma}\sqrt{\sigma^2+8a}}\right), \qquad \qquad \qquad \ &\!\!\!\!\!\sigma^2+8a>0, \\
				x^{\frac{1}{2}}e^{\frac{k}{\sigma^2}x+at}(c_1+c_2\ln x), \qquad \qquad \qquad \qquad \qquad \qquad \qquad \qquad &\!\!\!\!\!\sigma^2+8a=0, \\
				x^{\frac{1}{2}}e^{\frac{k}{\sigma^2}x+at}\left(c_1\cos(\frac{\sqrt{-\sigma^2-8a}}{2\sigma}\ln x)
				+c_2\sin(\frac{\sqrt{-\sigma^2-8a}}{2\sigma}\ln x)\right),
                \ &\!\!\!\!\! \sigma^2+8a<0. \!\end{array}\right.
\end{eqnarray}
The dynamic behavior of the pricing formula (\ref{inv8}) is shown in Figure 7 for different conditions of $\sigma^2+8a$. When $\sigma^2+8a>0$ (Figure 7 (a)), the price increases with
both $x$ and $t$. When $\sigma^2+8a=0$ (Figure 7 (b)), the price decreases with time $t$ but increases with $x$. When $\sigma^2+8a<0$ (Figure 7 (c)), the price initially rises and then falls with respect ot $x$, and first rises then stabilizes over time. These solutions provide insights into the pricing model's dynamics but do not account for boundary conditions. In practice, modifying the pricing model or changing its dependency or boundary conditions could make these solutions more applicable as pricing formulas.
\begin{figure}[H]
\centering
\subfigure[$\sigma^2 + 8 a > 0$, $a=1$.]{
		\begin{minipage}[b]{0.31\textwidth}
			\includegraphics[width=1\textwidth]{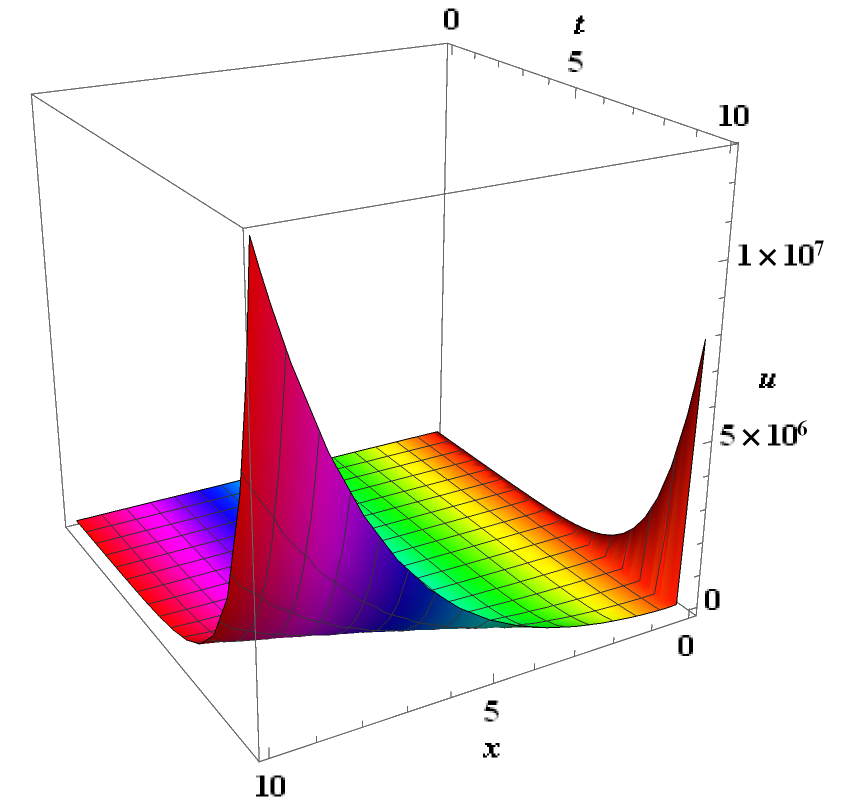}
		\end{minipage}
		\label{figsubd11}
	}
    \subfigure[$\sigma^2 + 8 a = 0$, $a=-0.5$.]{
    	\begin{minipage}[b]{0.32\textwidth}
   		\includegraphics[width=1\textwidth]{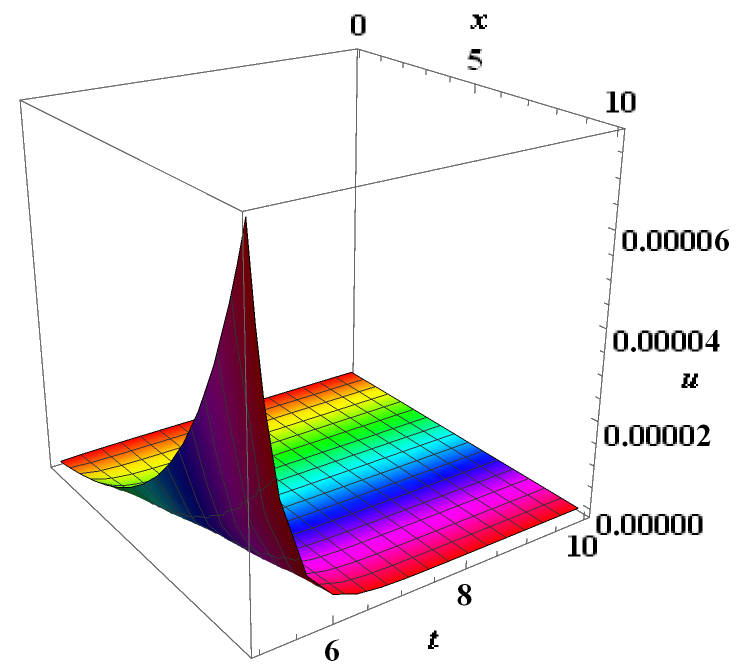}
    	\end{minipage}
	\label{figsubd12}
    }
    \subfigure[$\sigma^2 + 8 a < 0$, $a=-2.8$.]{
    	\begin{minipage}[b]{0.283\textwidth}
   		\includegraphics[width=1\textwidth]{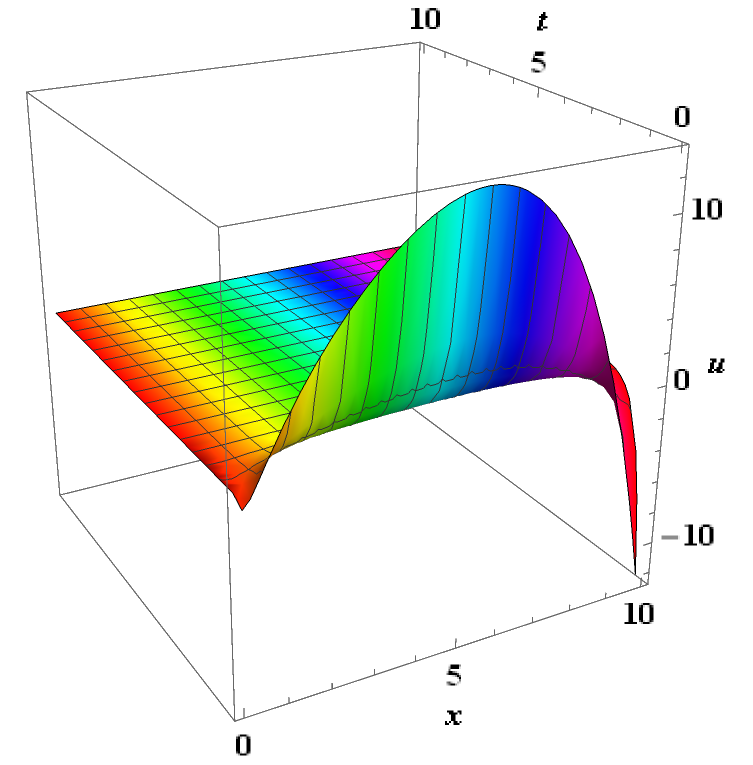}
    	\end{minipage}
	\label{figsubd13}
    }
	\caption{Three different dynamical structures of (\ref{inv8}) for $\sigma = 2$, $c_1 = 2$, $c_2 = 2$, $k = 1$.}\label{figd1}
\end{figure}
		
For the remaining two subalgebras, the one generated by $V_4$ has been previouly obtained, while the one generated by $V_6$ has no invariant solutions.

\medskip
Summarizing the discussions above, we arrive at a classification theorem for the invariant solutions of equation (\ref{eq1}).

\begin{theorem}
	An optimal system of invariant solutions of (\ref{eq1}) is constituted by solutions (\ref{c4}) and (\ref{inv5})-(\ref{inv8}), from which all other invariant solutions can be derived by suitable group transformations.
\end{theorem}

\section{Conclusion}\label{sec5}

In this paper, we performed a symmetry analysis for the Feynman-Kac formula associated with the GMR process, which is prevalent in the field of mathematical finance. We demonstrated that (\ref{eq1}) possesses a six-dimensional symmetry algebra plus an additional infinite-dimensional subalgebra under specific parameter conditions. By employing the representation of Lie algebra, we found an optimal system of one-dimensional subalgebras derived from the six-dimensional symmetry algebra. With the help of the optimal system we found, we performed symmetry reductions and constructed an optimal system of invariant solutions of (\ref{eq1}). We also offered interpretations of the financial behavior for several of these invariant solutions. The solutions we obtained could potentially be applied to model discounted project values or commodity prices within the finance market.

\section*{Acknowledgments}
The authors would like to thank the reviewers for their valuable comments and suggestions which greatly improved the work.









\medskip
Received August 2024;  revised October 2024; early access November 2024.
\medskip

\end{document}